\documentclass{amsart}
\usepackage{amsmath}
\usepackage{graphicx}
\usepackage{latexsym}

\theoremstyle{definition}

\newcommand{\CP}{{\mathbb{CP}}{}^{2}}

\begin{document}

\title[Handlebody argument for modifying achiral Lefschetz singularities]
{Handlebody argument for \\ modifying achiral Lefschetz singularities}

\author[R. \.{I}. Baykur]{R. \.{I}nan\c{c} Baykur}
\address{Department of Mathematics, Brandeis University, MA}
\email{baykur@brandeis.edu}

\date{06.17.2008.}


\maketitle

In \cite{B1} we gave a handlebody description of a broken Lefschetz fibration on $\CP$ as a counterexample to Gay and Kirby's conjecture on the necessity of negative Lefschetz singularities for generalized fibrations on arbitrary $4$-manifolds, and pointed out how this picture could be used to modify any given broken achiral Lefschetz fibration to a genuine broken Lefschetz fibration. Our general argument makes use of the following handlebody picture of a broken Lefschetz fibration over a disk, which can replace a regular neighborhood of a fiber with negative Lefschetz singularity:

\begin{figure}[ht] 
\begin{center}
\includegraphics[scale=0.85]{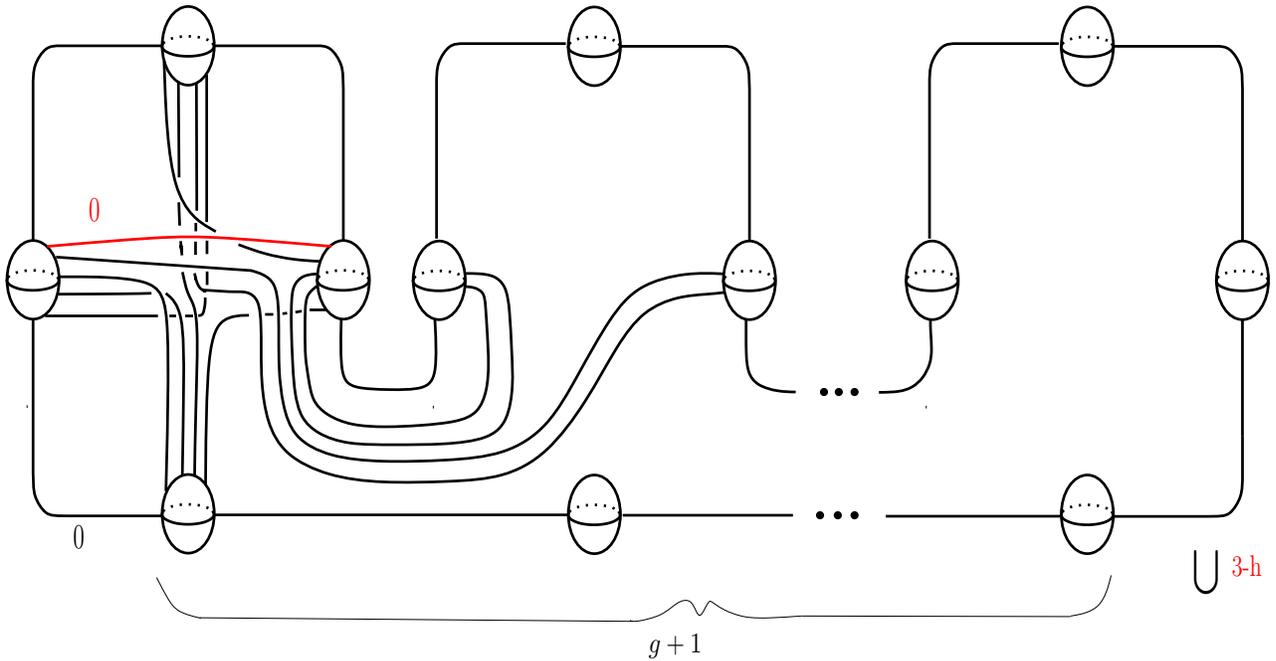}
\caption{\small The broken Lefschetz fibration over a disk to replace the neighborhood of a negative node.} \label{AchiralModification}
\end{center}
\end{figure}

\noindent Here the diagram is drawn from the `higher side', where there are three Lefschetz handles and a round $2$-handle attached to a fiber of genus $g+1$; so over the boundary of the base disk we have fibers of genus $g$. The $2$-handles corresponding to Lefschetz handles have fiber framing minus one, and the $2$-handle of the round $2$-handle given in red has fiber framing zero. The reader can turn to \cite{B0} for the conventions we use to depict broken Lefschetz fibrations in this note. 

\begin{figure}[ht] 
\begin{center}
\includegraphics[scale=0.9]{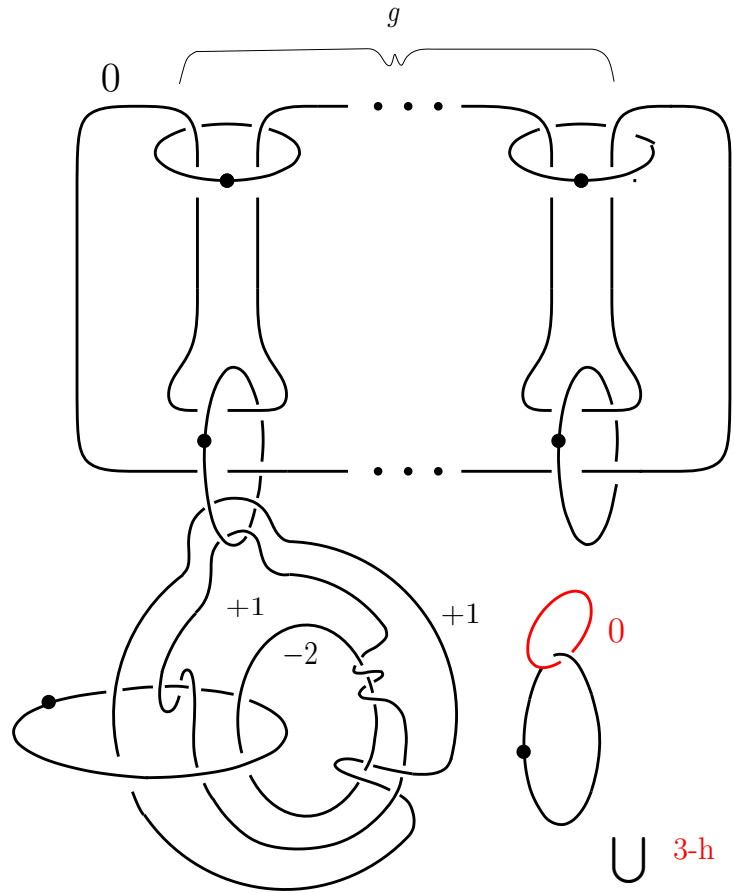}
\caption{\small Simplification of the handle diagram.} \label{AchiralSimplification}
\end{center}
\end{figure}

To show that both fibrations have the same total space we proceed as follows: Using the $0$-framed $2$-handle of the round $2$-handle we split the diagram as in Figure \ref{AchiralSimplification}, where we also switch to the dotted-circle notation to perform the rest of our handle calculus. We can now cancel the $2$-handle of the round $2$-handle against the $1$-handle it is linked with. Keeping the upper part of the diagram with $2g$ $1$-handles and the $0$-framed $2$-handle as it is, we will simplify the remaining part of the diagram where we have one $1$-handle, three $2$-handles, and a $3$-handle. We first slide the $(+1)$-framed $2$-handles over the $(-2)$-framed $2$-handle to separate them from the bottom left $1$-handle, and cancel this $1$-handle against this $(-2)$-framed $2$-handle. We then get two $(+1)$-framed unknots, linking once. One more handle slide separates an unknotted $2$-handle with framing $0$, and this $2$-handle can be canceled against the $3$-handle. The $(+1)$-framed $2$-handle we are left with still links with the $1$-handle contained in the upper part of the diagram, thus giving us a handlebody picture of a genus $g$ fibration over a disk and with one negative Lefschetz singularity.

Relying on the construction of Gay and Kirby in \cite{GK}, the above argument yields a handlebody proof of the existence of broken Lefschetz fibrations on arbitrary closed smooth oriented $4$-manifolds. (Another handlebody proof of this existence result was later given by Akbulut and Karakurt, where the achirality is avoided in a rather different way; see \cite{AK}.) The purpose of the current note is to reconstruct our picture locally on a twice punctured torus in the fiber and with a $3$-fold symmetry so as to provide a comparable picture with Lekili's achiral modification argument that uses singularity theory, given in Section 6 of his paper \cite{Lek}. One can certainly localize the handlebody picture we had in Figure \ref{AchiralModification} by throwing away the $0$-framed $2$-handle corresponding to the fiber and all the $1$-handles but the three $1$-handles the Lefschetz $2$-handles are linked with. However, in order to achieve the $3$-fold symmetry, we also need to rearrange the Lefschetz $2$-handles. This symmetric picture and many observations contained in the following paragraphs arose during the author's conversations with David Gay. 

\begin{figure}[ht] 
\begin{center}
\includegraphics[scale=0.7]{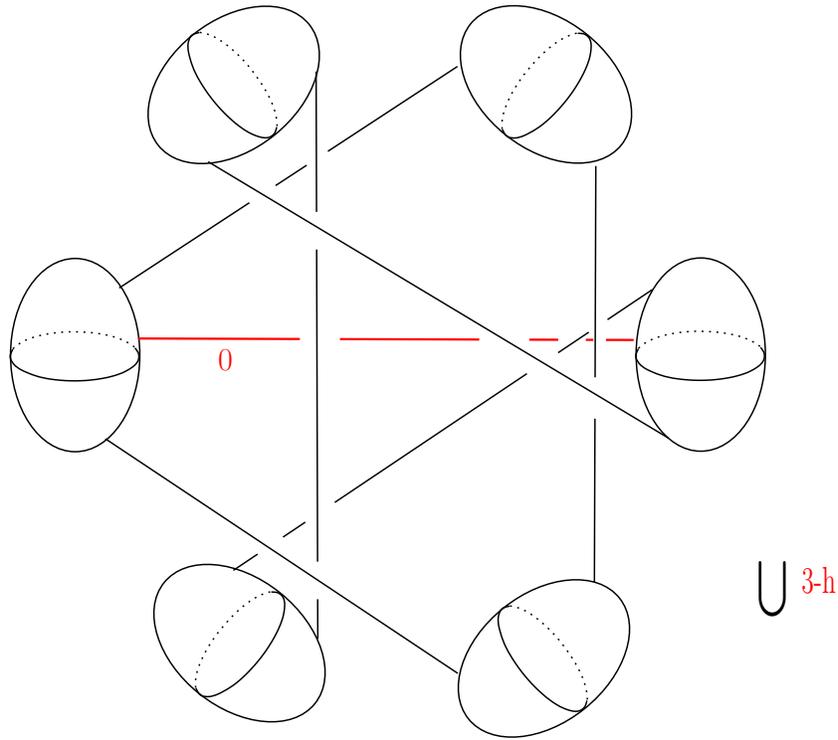}
\caption{\small The handlebody picture with a $3$-fold symmetry, which is obtained by rotating the diagram $2 \pi / 3$ degrees clockwise while shifting the order of the Lefschetz $2$-handles (towards the page) by one. All the $2$-handles but the $2$-handle of the round $2$-handle given in red have fiber framing minus one.} \label{HandlePicture}
\end{center}
\end{figure}

\begin{figure}[ht] 
\begin{center}
\includegraphics[scale=0.7]{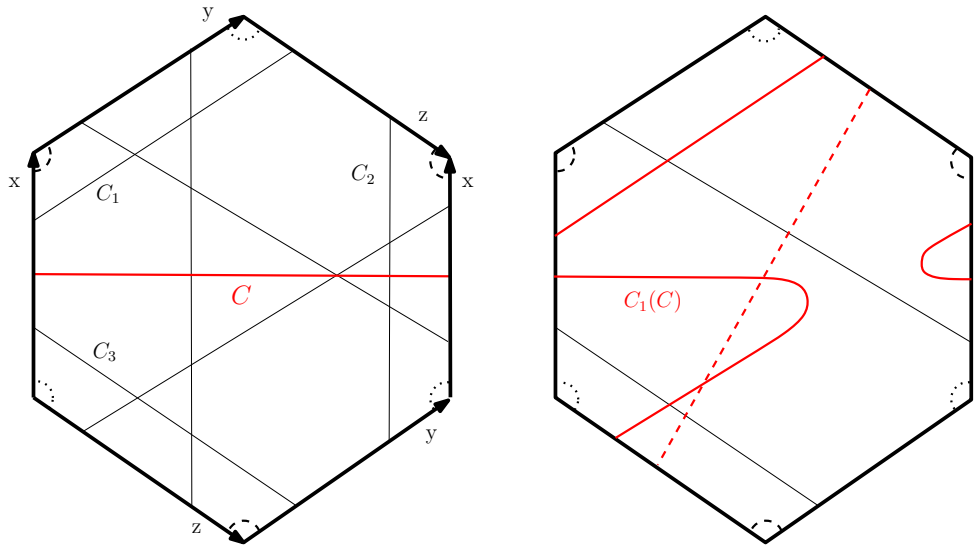}
\caption{\small On the left: The vanishing cycles on the fundamental domain of a twice punctured torus. On the right: The image of $C$ under the right-handed Dehn twist along $C_1$, represented by $C_1(C)$, can be isotoped on the fiber to the dashed red curve.} \label{Stabilization}
\end{center}
\end{figure}

In Figure \ref{HandlePicture} we see the twice punctured torus fiber, along with the vanishing cycles of three (positive) Lefschetz handles and a red curve which is the $2$-handle of the round $2$-handle. The $3$-handle completes the round $2$-handle. To show that this picture indeed gives a broken Lefschetz fibration over a disk one only needs to check that the attaching sphere of the $2$-handle of the round $2$-handle is sent to itself when the three Dehn twists prescribed by the Lefschetz handles are applied to it. Let us label the curves on the fiber as in Figure \ref{Stabilization}. After applying the first right-handed Dehn twist along $C_1$ one sees that the image of $C$ can be isotoped so that it is positioned with respect to the curve $C_2$ as $C$ was positioned with respect to $C_1$ in the first place (Figure \ref{Stabilization}). From the obvious $3$-fold symmetry one can conclude that after applying the right-handed Dehn twists along $C_1$, $C_2$ and $C_3$, the curve $C$ gets mapped onto itself (with the same orientation).

\begin{figure}[ht] 
\begin{center}
\includegraphics[scale=0.9]{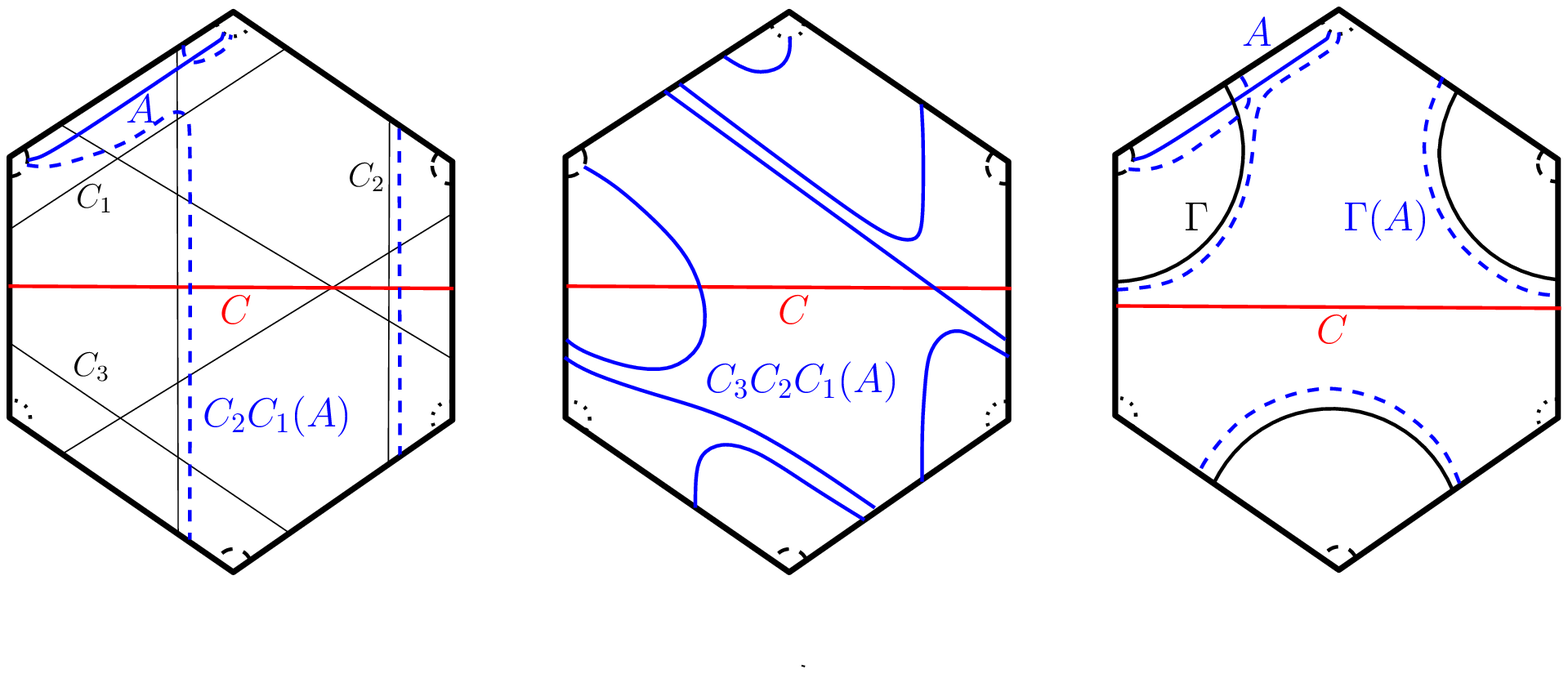}
\caption{\small Calculation of the reduced monodromy. } \label{ReducedMonodromy}
\end{center}
\end{figure}

It remains to verify two things: First is to see that the total space of this fibration (which in fact is the $4$-ball) is the same as that of a fibration with an annulus fiber and a single negative Lefschetz singularity attached along a separating curve on an annulus. The required calculus for this is similar to the one we have given above, and will be left to the reader. Secondly, we double check that the reduced monodromy of the fibration over the boundary of the base disk in the former fibration is equivalent to that of the latter. This allows us to interchange these pieces while matching the fibrations along the boundary, and therefore to extend the given fibrations on the rest. For this, recall that the mapping class group of the annulus is generated by a Dehn twist along a boundary parallel curve. So it suffices to understand the effect of the reduced monodromy on a simple arc $A$ that runs from one boundary component to another; see Figure \ref{ReducedMonodromy}. The curve $C_1$ does not intersect $A$, so the first nontrivially acting curve is $C_2$. The image of $A$ under the Dehn twist along $C_2$ is given on the first picture in Figure \ref{ReducedMonodromy} by the dashed blue curve. The blue curve $C_3 C_2 C_1 (A)$ in the second picture is the image of $A$ after all three Dehn twists are applied. This is where the round $2$-handle gets into action. We can slide the curve $C_3 C_2 C_1 (A)$ over the red curve (twice) and get the dashed blue curve given in the third picture in Figure \ref{ReducedMonodromy}. However, the third picture describes the effect of a \textit{left-handed} Dehn twist along the curve $\Gamma$, which is a boundary parallel curve on the annulus obtained after surgering the punctured torus along the red curve. Hence the two monodromies are the same. 

\begin{figure}[ht] 
\begin{center}
\includegraphics[scale=0.6]{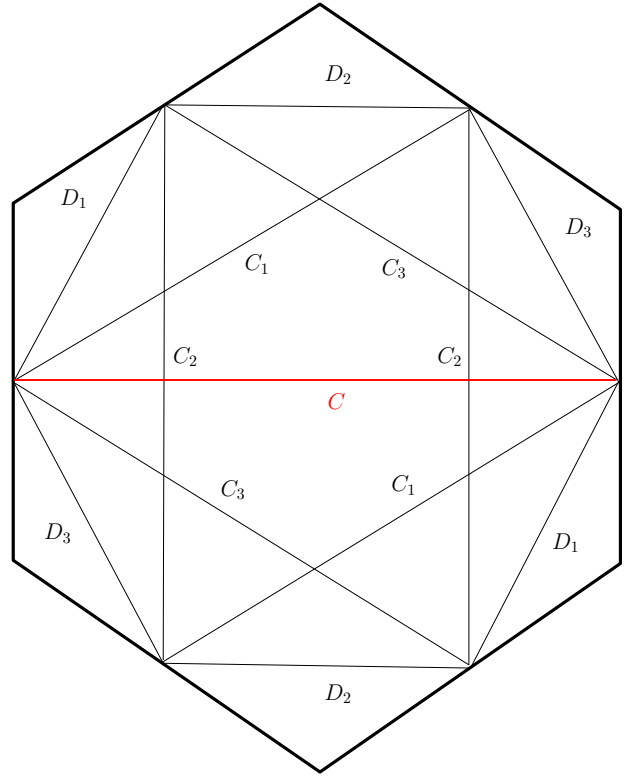}
\caption{\small ``Full symmetry.'' Vanishing cycles used in the modification around both positive and negative Lefschetz singularities are given.} \label{FullSymmetry}
\end{center}
\end{figure}

We would like to finish with a few observations. The Figure \ref{FullSymmetry} drawn in the most symmetric fashion presents a different choice of three vanishing cycles $D_1$, $D_2$ and $D_3$ on the twice punctured torus, and through similar arguments as above one can see that this picture stands for a broken Lefschetz fibration that can be used to replace a \textit{positive} Lefschetz singularity locally. This contains the nontrivial part of Perutz's example of a broken Lefschetz fibration (Example 1.3 in \cite{P}; also see Example 3.2 in \cite{B0}), and corresponds to the broken Lefschetz fibration that Lekili obtains after perturbing a positive Lefschetz singularity in his paper. One can then draw the curves $C_1$, $C_2$ and $C_3$ by joining the vertices of the hexagon formed by $D_1$, $D_2$, $D_3$ in the center as in Figure \ref{FullSymmetry}, and get the picture we had above (Figure \ref{Stabilization}, on the left) up to isotopy. From the very symmetry of the picture it is now easy to generalize our constructions to twice punctured $(4n+2)$-gons, for any $n \geq 1$, and thus to obtain various broken Lefschetz fibrations with higher fiber genera over disks. 

\begin{thebibliography}{99999}

\bibitem{AK} S. Akbulut and C. Karakurt, {\em  ``Every $4$-manifold is BLF'' }, J. of Gokova Geom. Topol. 2
(2008) 40-82. 

\bibitem{B0} R.\,I. Baykur, {\emph ``Topology of broken Lefschetz fibrations and near-symplectic four-manifolds'',} Pacific J. Math. 240 (2009), no. 2, 201--230; http://arxiv.org/abs/0801.0192.

\bibitem{B1} R.\,I. Baykur, {\em ``Existence of broken Lefschetz fibrations'',} Int. Math. Res. Not. IMRN 2008, Art. ID rnn 101, 15 pp. 

\bibitem{GK}D.\,T. Gay and R. Kirby, {\em ``Constructing Lefschetz-type fibrations on four-manifolds'',}  Geom. Topol. {\bf 11} (2007), 2075–-2115.

\bibitem{Lek} Y. Lekili, {\em ``Wrinkled Fibrations on near-symplectic manifolds'',} Geom. Topol. 13 (2009), 277-318.

\bibitem{P} T. Perutz, {\em ``Lagrangian matching invariants for fibred four-manifolds: I''}, Geom. Topol. {\bf 11} (2007), 759--828.
\end{thebibliography}
\end{document}